\documentclass{amsart}
\usepackage[english]{babel}

\usepackage{algorithm,algorithmic}
\usepackage{version}

\usepackage{ifpdf}
\ifpdf
  \usepackage[pdftex]{graphicx}
  \graphicspath{{./}}
  \DeclareGraphicsExtensions{.pdf}
  \usepackage[pdftex,hypertexnames=false,colorlinks=true,linkcolor=black,citecolor=black]{hyperref}
\else
  \usepackage[dvips]{graphicx}
  \graphicspath{{./}}
  \DeclareGraphicsExtensions{.eps}
  \usepackage[dvips,hypertexnames=false,colorlinks=true,linkcolor=black,citecolor=black]{hyperref}
\fi

\def\EuScript{\mathcal}
\def\mathscr{\EuScript}

\newcommand{\bbR}{\mathbb{R}}                               

\newcommand{\argmin}{\mathop{\arg\min}}                     

\newcommand{\espacea}[1]{\mathbb{#1}}                       
\newcommand{\tribu}[1]{\mathscr{#1}}                        
\newcommand{\omeg}{\Omega}                                  
\newcommand{\trib}{\tribu{A}}                               
\newcommand{\prbt}{\mathbb{P}}                              
\newcommand{\espe}{\mathbb{E}}                              

\newcommand{\va}[1]{\boldsymbol{\uppercase{#1}}}            

\newcommand{\vardelim}[1]{\left(#1\right)}                  

\newcommand{\esper}[2][]{\espe_{#1}\vardelim{#2}}           

\theoremstyle{remark}
\newtheorem{rem}{Remark}
\theoremstyle{plain}
\newtheorem{prop}{Proposition}

\begin{document}

\title{Decomposition of large-scale stochastic optimal control problems}

\author{Kengy Barty}
\address{Kengy Barty, EDF R\&D, 1, avenue du G\'{e}n\'{e}ral de Gaulle, 92141 Clamart Cedex, France}

\author{Pierre Carpentier}
\address{Pierre Carpentier, \'{E}cole Nationale Sup\'{e}rieure de Techniques Avanc\'{e}es (ENSTA), 32, boulevard Victor, 75015 Paris, France}

\author{Pierre Girardeau}
\address{Pierre Girardeau, Universit\'{e} Paris-Est, CERMICS, \'{E}cole des Ponts, Champs sur Marne, 77455 Marne la Vall\'{e}e Cedex 2, France, also with EDF R\&D and ENSTA}

\date{\today}

\begin{abstract}
	In this paper, we present an Uzawa-based heuristic that is adap\-ted to some type of stochastic optimal control problems. More precisely, we consider dynamical systems that can be divided into small-scale independent subsystems, though linked through a static almost sure coupling constraint at each time step. This type of problem is common in production/portfolio management where subsystems are, for instance, power units, and one has to supply a stochastic power demand at each time step. We outline the framework of our approach and present promising numerical results on a simplified power management problem.
\end{abstract}

\maketitle

\section{Introduction} \label{sec:Intro}

Stochastic optimal control is concerned about finding strategies to manage dynamical systems in an optimal way, with respect to some cost function. The particularity of such optimization problems is that the optimization variables we deal with are random variables. Indeed, the dynamical systems we consider are partially driven by some exogeneous noises and the objective function may also include such noises. Hence controls are random variables. Classical approaches such as Dynamic Programming and Stochastic Programming, that we briefly recall below, encounter difficulties when the system becomes large. The aim of this paper is to present a new heuristic to solve a class of such problems using price decomposition.

The problems we are studying are common in practice. For example, consider a physical system, say a set of numerous power units, that evolve depending on exogeneous noises (water inflows, failures) and on controls (production levels). At each time step, an observation on the system arises and a control has to be chosen on the basis of the available information, namely the past observations (non-anticipativity constraint). The objective is to minimize the sum of the units' production costs over a given discretized time horizon, while satisfying a global demand constraint at each time step. This decision process hence consists of finding optimal strategies, i.e. functions that map, at each time $t$, the available information to the optimal decision with respect to the production cost.

As far as we know, most methods that have been proposed to decompose large-scale stochastic optimal control problems are based on Stochastic Programming (see \cite{PrekopaStoProg,StochasticProgramming03}). This approach consists in representing the non-anticipativity constraints using a so-called scenario tree. Once discretized on such a structure, the problem is not stochastic anymore and various deterministic decomposition techniques have been used to solve it (see \cite{StochasticDecomposition96}). In this context, there are two main issues that are not easy to deal with. The first is concerned with the ``distance'' between the original problem and its deterministic reformulation \cite{Artstein91,HeitschRomischStrugarek05}, or how to draw a scenario tree in such a way that the solution of the discretized problem is an accurate estimate of the original one. In order to obtain some given accuracy, \cite{Shapiro06} shows that the growth of the number of leaves in the tree, hence the numerical complexity, has to be exponential with respect to the time horizon. The second issue is concerned with the way one can rebuild strategies from optimal commands obtained in the discretized problem \cite{Pennanen05}.

On the other hand, when dealing with a Markov Decision Process, methods based on Dynamic Programming (DP) (see \cite{Bellman57,BertsekasDP}) do provide a way to obtain strategies as feedback functions with respect to so-called state variables. Unfortunately, the well-known curse of dimensionality prevents us from using this approach straightforward on large-scale problems, because the computational burden increases exponentionally with the state dimension. Numerous approximations have been proposed to tackle the difficulty. For instance, a popular idea in the field of hydro-power management, introduced by Turgeon in \cite{Turgeon80}, consists of obtaining local strategies as a function of the local stock and the aggregated complementary stock. Another idea, namely Approximate Dynamic Programming (ADP), looks for the value functions (solutions of the DP equation) within a finite-dimensional space (see \cite{BellmanDreyfus59} or \cite[\S6.5]{BertsekasTsitsiklis96}). To be practically efficient, such an approach requires some a priori information about the problem, in order to define a well suited functional subspace. Indeed, there is no systematic means to choose the basis functions and several choices have been proposed \cite{TsitsiklisVanRoy96,deFariasVanRoy03}.

When dealing with large-scale optimization problems, the decomposition/co\-ordi\-nation approach aims at finding a solution to the original problem by iteratively solving several smaller-dimensional subproblems. In the deterministic case, several types of decomposition have been proposed (e.g. by prices or by quantities) and unified in a general framework using the Auxiliary Problem Principle in \cite{Cohen80}. In the open-loop stochastic case, i.e. when controls do not rely on any observation, \cite{CohenCulioli90} proposed to take advantage of both decomposition techniques and stochastic gradient algorithms. These techniques have been extended in the closed-loop stochastic case by \cite{BartyRoyStrugarek05}, but so far they fail to provide decomposed state dependent strategies in the Markovian case. This is because a subproblem's optimal strategy depends on the state of the whole system, not only on the local state. In other words, decomposition approaches are meant to decompose the control space, namely the range of the strategy, but the numerical complexity of the problems we consider here also arises because of the dimensionality of the state space, that is to say the domain of the strategy.

We here propose a way to use price decomposition within the closed-loop stochastic case. The coupling constraints, namely the constraints preventing the problem from being naturally decomposed, are dualized using a Lagrange multiplier (price). At each iteration, the price decomposition algorithm solves each subproblem using the current prices, then uses the solutions to update the prices. In the stochastic context, prices are a random process whose dynamics are not available, so the subproblems do not in general fall into the Markovian setting. However, on a specific instance of this problem, \cite{TheseStrugarek} has exhibited a dynamics for the optimal multipliers, and he has shown that these dynamics were independent of the decision variables. Hence it was possible to come down to the Markovian framework and to use DP to solve the subproblems in this case. Following this idea, we propose to choose a parameterized dynamics for these multipliers so that solving subproblems using DP becomes possible. The update is then performed using a sampling/regression technique.

This paper is organized as follows. In \S\ref{sec:MF} we describe the general type of problems we are concerned with in this paper. Then, in \S\ref{sec:DP}, we recall the Dynamic Programming equation and highlight the difficulties induced when considering large-scale problems. In \S\ref{sec:PD}, we present the classical price decomposition approach in Hilbert spaces and the difficulties encountered when dealing with stochastic optimal control problems. Based on these ingredients, we present in \S\ref{sec:Point} a heuristic allowing us to solve subproblems using DP. We finally validate this approach on a simplified power management problem in \S\ref{sec:NumExp}.

\section{Mathematical framework} \label{sec:MF}

Throughout this paper the random variables, defined over a probability space $\left(\omeg, \trib, \prbt\right)$, will be denoted using bold letters (e.g. $\va{w} \in L^2\left(\omeg, \trib, \prbt, \espacea{W}\right)$) whereas their realizations will be denoted using normal letters (e.g. $w \in \espacea{W}$).

In this paper we consider a finite horizon stochastic optimal control problem, where $T$ denotes the time horizon. Three types of random variables are involved in the problem, namely a state, a control, and a noise. The state $\va{x}_t \in L^2\left(\omeg, \trib, \prbt, \bbR^n\right)$ evolves with respect to dynamics depending on the control $\va{u}_t \in L^2\left(\omeg, \trib, \prbt, \bbR^m\right)$ and on some exogeneous noise $\va{\xi}_t \in L^2\left(\omeg, \trib, \prbt, \bbR^p\right)$. Unlike deterministic optimal control problems, in the stochastic case the time ``direction'' is of particular importance. In order to fulfill the causality principle, the control at a given time step $t$ only depends on the observation of noises prior to $t$. Moreover, we assume that the observation available at time $t$ consists of all past noises. In order to mathematically represent such an information structure, we denote by $\trib_t$ the $\sigma$-field generated at time $t$ by past noises $\left(\va{\xi}_1, \dots, \va{\xi}_t\right)$, so that the control $\va{u}_t$ at time step $t$ has to be measurable with respect to $\trib_t$. These last constraints will be called the non-anticipativity constraints.

The global system consists of $N$ units, whose dynamics and cost functions are mutually independent. More precisely, the state $\va{x}_t$ (respectively the control $\va{u}_t$) of the global system writes $\big(\va{x}_t^1, \dots, \va{x}_t^N\big)$ with $\va{x}_t^i \in L^2\left(\omeg, \trib, \prbt, \bbR^{n_i}\right)$ (respectively $\big(\va{u}_t^1, \dots, \va{u}_t^N\big)$ with $\va{u}_t^i \in L^2\left(\omeg, \trib, \prbt, \bbR^{m_i}\right)$) and $n=\sum_{i=1}^N n_i$ and $m=\sum_{i=1}^N m_i$, so that the global dynamics $\va{x}_{t+1} = f_t\left(\va{x}_t, \va{u}_t, \va{\xi}_{t+1}\right)$ can be written independently unit by unit: $\va{x}_{t+1}^i = f_t^i\left(\va{x}_t^i, \va{u}_t^i, \va{\xi}_{t+1}\right)$, $i=1, \dots, N$. In the same way, the global cost $L_t\left(\va{x}_t, \va{u}_t, \va{\xi}_{t+1}\right)$ is equal to the sum of the local unit costs $L_t^i\left(\va{x}_t^i, \va{u}_t^i, \va{\xi}_{t+1}\right)$, $i=1, \dots, N$. At the end of the time period, each unit $i$ leads a cost $K^i$ that only depends on the final state $\va{x}_T^i$ of the unit.

For now, the global problem can be stated independently unit by unit. The coupling between the units arises from a set of static $\bbR^d$-valued constraints, the constraint at time step $t$ reading $\sum_{i=1}^N g_t^i\left(\va{x}_t^i, \va{u}_t^i\right) = 0$ (see remark \ref{rem:Coupling} for extensions to more enhanced relations). We suppose that all functions $f_t^i$, $L_t^i$ and $g_t^i$ are at least Borel measurable.

The initial state $\va{x}_0$ is assumed to be known. Denoting $\left(\va{x}_1, \dots, \va{x}_T\right)$ by $\va{x}$ and $\left(\va{u}_0, \dots, \va{u}_{T-1}\right)$ by $\va{u}$, the problem we wish to solve is:
\begin{subequations} \label{eqn:P}
	\begin{align}
		\min_{\va{x},\va{u}}\qquad & \esper{\sum_{t=0}^{T-1} \sum_{i=1}^N L_t^i\left(\va{x}_t^i, \va{u}_t^i, \va{\xi}_{t+1}\right) + \sum_{i=1}^N K^i\left(\va{x}_T^i\right)} \\
		\text{s.t.}\qquad & \va{x}_{t+1}^i = f_t^i\left(\va{x}_t^i, \va{u}_t^i, \va{\xi}_{t+1}\right), \qquad \forall t=0, \dots, T-1, \forall i=1, \dots, N, \label{eqn:P-Dyn} \\
		& \va{x}_0^i = x^i, \qquad \forall i=1, \dots, N, \\
		& \sum_{i=1}^N g_t^i\left(\va{x}_t^i, \va{u}_t^i\right) = 0, \qquad \forall t=0, \dots, T-1, \label{eqn:P-Dem} \\
		& \va{u}_t \text{ is } \trib_t \text{-measurable}, \qquad \forall t=0, \dots, T-1, \label{eqn:P-Meas} \\
		& \underline{x}_t \leq \va{x}_t \leq \overline{x}_t, \qquad \forall t=1, \dots, T, \\
		& \underline{u}_t \leq \va{u}_t \leq \overline{u}_t, \qquad \forall t=0, \dots, T-1.
	\end{align}
\end{subequations}

\medskip

There are three types of coupling in Problem \eqref{eqn:P}:
\begin{itemize}
	\item The first one comes from the state dynamics \eqref{eqn:P-Dyn} that induce a temporal coupling, subsystem by subsystem.
	\item The second one arises from the static constraints \eqref{eqn:P-Dem} that link together all the subsystems at each time step $t$.
	\item The third type of coupling comes from the non-anticipativity constraints \eqref{eqn:P-Meas}, which link together controls relying on the same noise history. If two realizations of the noise process are identical up to time $t$, then the same control has to be applied at time $t$ on both realizations.
\end{itemize}

We ultimately suppose that noises $\va{\xi}_t$ are independent (white noise). We are thus in the Markovian case and it is well known that the optimal control, which is a priori a function of all the past noises, only depends on the current state \cite{BertsekasDP}.

\begin{rem}[White noise assumption]
If the noises $\va{\xi}_t$ are not independent\footnote{and also the noises $\va{d}_t$ introduced in remark \ref{rem:Coupling}} but still have known dynamics, one can always include the necessary noise history in the state to come back to the Markovian case. Unfortunately, this usually leads to a higher state dimension, and hence a higher numerical complexity in the DP framework, as will be explained in \S\ref{sec:DP}.
\end{rem}

\begin{rem}[Coupling constraints involving noises] \label{rem:Coupling}
It is possible to replace the static coupling constraint $\sum_{i=1}^N g_t^i\left(\va{x}_t^i, \va{u}_t^i\right) = 0$ by $\sum_{i=1}^N g_t^i\left(\va{x}_t^i, \va{u}_t^i\right) = \va{d}_t$, where $\va{d}_t$ is a random variable representing for instance a global demand. However, expressions are then harder to write: in the Markovian case, i.e. when the $\va{d}_t$'s are independent one from another, $\va{d}_t$ is observed before choosing the control at time $t$, so optimal controls must depend on both the state $\va{x}_t$ and the noise $\va{d}_t$.
\end{rem}

\section{Stochastic Dynamic Programming} \label{sec:DP}

In order to solve stochastic optimal control problems in the Markovian framework, Bellman proposed in \cite{Bellman57} the Dynamic Programming (DP) method. It consists of introducing value functions $V_t: \bbR^n \rightarrow \overline{\bbR}$ that represent the expected optimal cost when starting from state $x$ at time $t$. In the case of Problem \eqref{eqn:P}, it reads:
\begin{equation} \label{eqn:DefValueFunction}
	V_t\left(x\right) = \min_{\va{x},\va{u}} \displaystyle\esper{\left.\sum_{s=t}^{T-1} \sum_{i=1}^N L_s^i\left(\va{x}_s^i, \va{u}_s^i, \va{\xi}_{s+1}\right) + \sum_{i=1}^N K^i\left(\va{x}_T^i\right) \right\vert \va{x}_t = x},
\end{equation}
subject to the same constraints as in \eqref{eqn:P} and using the convention that if the optimization problem \eqref{eqn:DefValueFunction} is not feasible, then $V_t\left(x\right) = +\infty$. The value functions are usually computed in a recursive manner using the DP equation:
\begin{subequations} \label{eqn:DP}
	\begin{align}
		V_T\left(x\right)\quad &= \sum_{i=1}^N K^i\left(x^i\right), \label{eqn:DP-T} \\
	\intertext{and, for $t=1, \dots, T-1$:}
		V_t\left(x\right)\quad &= \min_{u \in \left[\underline{u}_t, \overline{u}_t\right]} \esper{\sum_{i=1}^N L_t^i\left(x^i, u^i, \va{\xi}_{t+1}\right) + V_{t+1}\left(f_t\left(x, u, \va{\xi}_{t+1}\right)\right)}, \label{eqn:DP-t} \\
		&\qquad\quad \text{s.t.} \quad \sum_{i=1}^N g_t^i\left(x^i, u^i\right) = 0. \label{eqn:DP-Coup}
	\end{align}
\end{subequations}

Unlike Stochastic Programming methods, a major advantage of DP is that it provides the control $\va{u}_t$ as a feedback function on the state variable $\va{x}_t$:
\begin{equation*}
	\va{u}_t = \Phi_t\left(\va{x}_t\right).
\end{equation*}

Except on very simple examples, Equation \eqref{eqn:DP} cannot be solved analytically, and many numerical methods have been proposed. A common practice is to discretize the state space and estimate the expectations using Monte Carlo sampling. Unfortunately, as was mentioned in \S\ref{sec:Intro}, we are facing the curse of dimensionality: the complexity of DP grows exponentially with respect to the state space dimension.

Moreover, Equation \eqref{eqn:DP} is not decomposable in the sense that it cannot be replaced by the solving of $N$ DP equations depending only on the local state $x^i$. Indeed, even if $V_T$ is a sum of functions depending on the local state $x^i$ as in \eqref{eqn:DP-T}, this additive property does not hold for the preceeding time steps because of the coupling constraint \eqref{eqn:DP-Coup}. Hence, looking for the value function as a sum of functions depending only on the local state would lead to suboptimal strategies. In other words, the local state of a subsystem is not sufficient to take the optimal local decision; some global information about the system is necessary.

Nonetheless, DP remains a seductive approach for small-scale problems since it provides a way to obtain feedback functions. Based on a decomposition scheme presented in \S\ref{sec:PD}, we will describe in \S\ref{sec:Point} a heuristic approach in which Problem \eqref{eqn:P} is decomposed into small-scale subproblems that we solve using DP.

\section{Price decomposition} \label{sec:PD}

Let us recall some results about the classical Uzawa algorithm \cite{ArrowHurwiczUzawa58}, which aims at iteratively getting round the static coupling constraint \eqref{eqn:P-Dem}. When the cost function is additive, this algorithm is also referred to as the price decomposition approach (see \cite{Cohen80} for further details). Let us first introduce the Lagrangian of problem~\eqref{eqn:P}:
\begin{equation*}
	\mathcal{L}\left(\va{x}, \va{u}, \va{\lambda}\right) \! = \!
\esper{\sum_{t=0}^{T-1} \sum_{i=1}^N \left(L_t^i\left(\va{x}_t^i, \va{u}_t^i, \va{\xi}_{t+1}\right) + \va{\lambda}_t^\top g_t^i\left(\va{x}_t^i, \va{u}_t^i\right)\right) + \sum_{i=1}^N K^i\left(\va{x}_T^i\right)} \! ,
\end{equation*}
with $\va{\lambda}_t \in L^2\left(\omeg, \trib, \prbt, \bbR^d\right)$ the Lagrange multiplier associated to the coupling constraint \eqref{eqn:P-Dem} and $\va{\lambda} = \left(\va{\lambda}_0, \dots, \va{\lambda}_{T-1}\right)$. When the Lagrangian has a saddle point, we know from classical duality theory in optimization \cite{EkelandTemam92} that Problem \eqref{eqn:P} is equivalent to:
\begin{subequations} \label{eqn:Duality2}
	\begin{align}
		\max_{\va{\lambda}} \min_{\va{x},\va{u}}\quad &\mathcal{L}\left(\va{x}, \va{u}, \va{\lambda}\right) \\
		\text{s.t.}\quad & \va{x}_{t+1}^i = f_t^i\left(\va{x}_t^i, \va{u}_t^i, \va{\xi}_{t+1}\right), \qquad \forall t=0, \dots, T-1, \forall i=1, \dots, N, \label{eqn:Duality2-Dyn}\\
		& \va{u}_t \text{ is } \trib_t \text{-measurable}, \qquad \forall t=0, \dots, T-1, \label{eqn:Duality2-Meas}\\
		& \underline{x}_t \leq \va{x}_t \leq \overline{x}_t, \qquad \forall t=1, \dots, T, \label{eqn:Duality2-Bounds1} \\
		& \underline{u}_t \leq \va{u}_t \leq \overline{u}_t, \qquad \forall t=0, \dots, T-1. \label{eqn:Duality2-Bounds2}
	\end{align}
\end{subequations}

Recall that the Lagrange multiplier $\va{\lambda}_t$ can be interpreted as the marginal price one should pay for satisfying the coupling constraint \eqref{eqn:P-Dem}. Because of the $\trib_t$-measurability of the variables involved in this constraint and of the properties of conditional expectation, it is easy to see that we can always choose $\va{\lambda}_t$ to be $\trib_t$-measurable.

Let us introduce the dual function $\psi\left(\va{\lambda}\right) := \min_{\va{x},\va{u}} \mathcal{L}\left(\va{x}, \va{u}, \va{\lambda}\right)$ subject to constraints \eqref{eqn:Duality2-Dyn}, \eqref{eqn:Duality2-Meas}, \eqref{eqn:Duality2-Bounds1} and \eqref{eqn:Duality2-Bounds2}. The key point of the price decomposition algorithm is that computing $\psi\left(\va{\lambda}\right)$ is much easier than solving the original Problem \eqref{eqn:P}. Indeed, one can write:
\begin{align*}
	\psi\left(\va{\lambda}\right) &= \min_{\va{x},\va{u}} \esper{\sum_{t=0}^{T-1} \sum_{i=1}^N \left(L_t^i\left(\va{x}_t^i, \va{u}_t^i, \va{\xi}_{t+1}\right) + \va{\lambda}_t^\top g_t^i\left(\va{x}_t^i, \va{u}_t^i\right)\right) + \sum_{i=1}^N K^i\left(\va{x}_T^i\right)}, \\
	&= \sum_{i=1}^N \min_{\va{x}^i,\va{u}^i} \esper{\sum_{t=0}^{T-1} \left(L_t^i\left(\va{x}_t^i, \va{u}_t^i, \va{\xi}_{t+1}\right) + \va{\lambda}_t^\top g_t^i\left(\va{x}_t^i, \va{u}_t^i\right)\right) + K^i\left(\va{x}_T^i\right)},
\end{align*}
so that we replace the solving of an optimization problem with variables $\left(\va{x}, \va{u}\right)$ by the solving of $N$ subproblems with variables $\left(\va{x}^i, \va{u}^i\right)$.

Given $\va{\lambda}^k$, an iteration of the price decomposition algorithm first solves the $N$ subproblems:
\begin{subequations} \label{eqn:Local}
	\begin{align}
		\min_{\va{x}^i,\va{u}^i}\quad &\esper{\sum_{t=0}^{T-1} \left(L_t^i\left(\va{x}_t^i, \va{u}_t^i, \va{\xi}_{t+1}\right) + {\va{\lambda}_t^k}^\top g_t^i\left(\va{x}_t^i, \va{u}_t^i\right)\right) + K^i\left(\va{x}_T^i\right)} \label{eqn:Local1} \\
		\text{s.t.}\quad & \va{x}_{t+1}^i = f_t^i\left(\va{x}_t^i, \va{u}_t^i, \va{\xi}_{t+1}\right), \qquad \forall t=0, \dots,  T-1, \label{eqn:Local2} \\
		& \va{u}_t^i \text{ is } \trib_t \text{-measurable}, \qquad \forall t=0, \dots, T-1, \label{eqn:Local3} \\
		& \underline{x}_t^i \leq \va{x}_t^i \leq \overline{x}_t^i, \qquad \forall t=1, \dots, T, \\
		& \underline{u}_t^i \leq \va{u}_t^i \leq \overline{u}_t^i,\qquad \forall t=0, \dots, T-1.
	\end{align}
\end{subequations}
The Lagrange multiplier $\va{\lambda}^k$ is then updated using a gradient-like algorithm. Under standard assumptions\footnote{See \cite{Danskin67} for results on the differentiability of the dual fonction $\psi$.}, the gradient of $\psi$ is:
\begin{equation*}
	\nabla_{\va{\lambda}_t}\psi\left(\va{\lambda}^k\right) = \sum_{i=1}^N g_t^i\left(\va{x}_t^{i,k+1}, \va{u}_t^{i,k+1}\right),
\end{equation*}
where $\va{x}_t^{i,k+1}$ and $\va{u}_t^{i,k+1}$ are the solutions of Problem \eqref{eqn:Local}.

At first sight, Problem \eqref{eqn:Local} looks like a stochastic optimal control problem with control $\va{u}_t^i$ and state $\va{x}_t^i$, the solution of which would be a local feedback on $\va{x}_t^i$. This contradicts the fact that the solution of Problem \eqref{eqn:P} is a feedback function on the whole state $\big(\va{x}_t^1, \dots, \va{x}_t^N\big)$. In order to understand where this contradiction comes from, one has to highlight the role of $\va{\lambda}$ in Problem \eqref{eqn:Local}.

\section{Dual Approximate Dynamic Programming} \label{sec:Point}

Let us take a closer look at Problem \eqref{eqn:Local}. First suppose that $\va{\lambda}$ is a white noise process. Then Problem \eqref{eqn:Local} lies in the Markovian framework with state $\va{x}_t^i$ and noise $\left(\va{\xi}_t, \va{\lambda}_t\right)$. The optimal control $\va{u}_t^i$ depends only on the local state $\va{x}_t^i$ and one can apply stochastic dynamic programming to solve this small-scale optimal control problem. Unfortunately, we do not know anything about the time correlations of the price process $\va{\lambda}$\ldots

Let us now consider the general case. Defining $\left(\va{x}_t^i, \va{\lambda}_1, \dots, \va{\lambda}_t\right)$ as the state at time $t$, Problem \eqref{eqn:Local} falls in the Markovian setting. In particular, the optimal control $\va{u}_t^i$ is $\left(\va{x}_t^i, \va{\lambda}_1, \dots, \va{\lambda}_t\right)$-measurable. However DP in this context proves numerically intractable because the state dimension increases with respect to time.

Consider now an intermediate case, and suppose that the dual variable $\va{\lambda}$ has a short memory dynamics, for instance that $\va{\lambda}_{t+1}$ only depends on $\va{\lambda}_t$ and $\va{\xi}_{t+1}$:
\begin{equation}
	\va{\lambda}_{t+1} = h_t\left(\va{\lambda}_t, \va{\xi}_{t+1}\right).
\end{equation}
Using $\left(\va{x}_t^i, \va{\lambda}_t\right)$ as the state variable at time $t$, Problem \eqref{eqn:Local} falls in the Markovian setting. The state dimension does not increase with respect to time anymore and is hopefully small so that Problem \eqref{eqn:Local} can be solved using DP.

In a very specific instance of Problem \eqref{eqn:P}, namely:
\begin{equation} \label{eqn:PStrugarek}
	\begin{array}{rl}
		\displaystyle\min_{\va{x},\va{u}} & \displaystyle\esper{\sum_{t=0}^{T-1} \sum_{i=1}^N \frac{c_i}{2} \left(\va{u}_t^i\right)^2 + \sum_{i=1}^N \frac{\gamma_i}{2} \left(\va{x}_T^i-\va{x}_0^i\right)^2} \\
		\text{s.t.} & \va{x}_{t+1}^i = \va{x}_t^i - \va{u}_t^i + \va{a}_{t+1}^i, \qquad \forall t=0, \dots, T-1, \\
		& \displaystyle\sum_{i=1}^n \va{u}_t^i = \va{d}_t, \qquad \forall t=0, \dots, T-1, \\
		& \va{u}_t \text{ is } \trib_t \text{-measurable}, \qquad \forall t=0, \dots, T-1,
	\end{array}
\end{equation}
with $\va{\xi}_t = \big(\va{d}_t, \va{a}^1_t, \dots, \va{a}^N_t\big)$, \cite{TheseStrugarek} has brought to light such an intermediate case. Here the dimension of the state $\va{x}_t^i$ (respectively of the control $\va{u}_t^i$) in the subsystem $i$ is $n_i=1$ (respectively $m_i=1$), for $i=1, \dots, N$. The result is the following.
\begin{prop} \label{prop:Strugarek}
	If $\big(\va{d}_t, \va{a}_t^1, \dots, \va{a}_t^N\big)_{t=0, \dots, T}$ is a white noise process and if there exists $\alpha \in \bbR^+$ such that $\gamma_i = \alpha c_i, \forall i=1, \dots, N$, then the optimal Lagrange multipliers satisfy:
	\begin{equation*}
		\begin{array}{rl}
			\va{\lambda}_{t+1} &= \displaystyle\va{\lambda}_t + \frac{1}{\sum_{i=1}^N \frac{1}{c_i}} \Big(\va{d}_{t+1} \left(1+\alpha\right) - \va{d}_t - \alpha \esper{\va{d}_{t+1}} \\
			&\qquad - \alpha \left(\va{a}_{t+1}-\esper{\va{a}_{t+1}}\right)\Big), \\
			\va{\lambda}_0 &= \displaystyle\frac{1}{\sum_{i=1}^N \frac{1}{c_i}} \left(\va{d}_0 \left(1-\alpha\right) - \alpha \sum_{s=1}^T \esper{\va{a}_s} - \alpha \sum_{s=1}^{T-1} \esper{\va{d}_s}\right).
		\end{array}
	\end{equation*}
\end{prop}
Using such a dynamics for the multipliers, it is straightforward to show that Problem \eqref{eqn:PStrugarek} splits into $N$ independent optimization subproblems. Taking the state variable as $\left(\va{x}_t^i, \va{\lambda}_t, \va{d}_t\right)$, the $i$-th subproblem can be solved using DP in dimension 3. In summary, we have replaced one $N$-dimensional problem by $N$ $3$-dimensional problems.

Note that the proportionality assumption on the cost coefficients in proposition \ref{prop:Strugarek} is rather unnatural. Nevertheless, it shows that, in some cases, there exist dynamics for the Lagrange multipliers that is independent of the decision variables.

To deal with more general cases, we propose to approximate the dual process $\va{\lambda}$ by some parameterized short-memory process. That is, we try to identify the multipliers that are the closest to the optimal ones within a constrained subspace of stochastic processes. This approach is similar to that employed in the Approximate Dynamic Programming (ADP) method. Since it concerns dual variables rather than DP value functions, we refer to this approach as Dual Approximate Dynamic Programming.

The performance of such an approach highly depends on the choice of the subspace of stochastic processes in which we force the multipliers to lie. However, a major difference with ADP techniques is that approximating the dual variables may lead to violations of the coupling constraints. The larger the chosen subspace, the less the coupling constraints will be violated. Moreover, prior information on the problem may be useful to devise a suitable dynamics.

Let us now present the implementation of DADP. We constrain dual variables to satisfy:
\begin{equation} \label{eqn:DualDynamics}
	\va{\lambda}_{t+1} = h_{\alpha_t}\left(\va{\lambda}_t, \va{\xi}_{t+1}\right),
\end{equation}
where $h_{\alpha_t}$ is an a priori chosen function parameterized by $\alpha_t \in \bbR^q$. We denote by $\mathcal{S}$ the set of all random processes that verify Equation \eqref{eqn:DualDynamics} for some real vector $\alpha=\left(\alpha_1, \dots, \alpha_{T-1}\right)$. Given a vector $\alpha^{k}$ of coefficients, the first step of DADP is to solve the $N$ subproblems \eqref{eqn:Local} using DP with state $\left(\va{x}_t^i, \va{\lambda}_t\right)$. In order to update the Lagrange multipliers, we draw $s$ trajectory samples of the noise $\va{\xi}$ and integrate the dynamics \eqref{eqn:Local2} and \eqref{eqn:DualDynamics} using the optimal feedback laws, thus obtaining $s$ trajectory samples of $\va{x}^k$, $\va{u}^k$ and $\va{\lambda}^k$. We then perform a gradient step on $\va{\lambda}$ sample by sample:
\begin{equation*}
	\va{\lambda}_t^{k+\frac{1}{2}, \sigma} = \va{\lambda}_t^{k, \sigma} + \rho_{t} \sum_{i=1}^N g_t^i\left(\va{x}_t^{i,k,\sigma}, \va{u}_t^{i,k,\sigma}\right), \qquad \forall \sigma=1, \dots, s,
\end{equation*}
with $\rho_{t}$ being well-chosen real values. Finally, we apply a regression operator $\mathcal{R}^s$ on the samples $\va{\lambda}^{k+\frac{1}{2}}$ in order to obtain a stochastic process $\va{\lambda}^{k+1}$ lying in $\mathcal{S}$:
\begin{align*}
	\mathcal{R}^s\left(\va{\lambda}^{k+\frac{1}{2}}\right) = &\argmin_{\alpha_0, \dots, \alpha_{T-1}} \sum_{t=0}^{T-1} \sum_{\sigma=1}^s \left\Vert\va{\mu}_t^\sigma - \va{\lambda}_t^{k+\frac{1}{2},\sigma} \right\Vert_{\bbR^d}^2 \\
	&\text{s.t.} \quad \va{\mu}_{t+1}^\sigma = h_{\alpha_t}\left(\va{\mu}_t^\sigma, \va{\xi}_{t+1}^\sigma\right).
\end{align*}

This heuristic is outlined in algorithm \ref*{alg:Heuristic}.
\begin{algorithm}
\caption{Dual Approximate Dynamic Programming}
\label{alg:Heuristic}
\begin{algorithmic}
	\REQUIRE $\varepsilon > 0$, $\gamma > 0$, a shape $h_\alpha$ for the prices dynamics, $\alpha^0$.
	\REPEAT
		\STATE $k \leftarrow k+1$
		\FOR{$i=1$ to $N$}
			\STATE Solve $i$-th subproblem by DP using parameters $\alpha^k$ for the price dynamics, and obtain $\va{x}^{i,k}$ and $\va{u}^{i,k}$. Both implicitly depend on $\alpha^k$.
		\ENDFOR
		\STATE Update paramaters $\alpha^k$:
		\begin{align*}
			\alpha^{k+1} &= \mathcal{R}^s\left(\left(\va{\lambda}_t^k + \rho_{t} \sum_{i=1}^N g_t^i\left(\va{x}_t^{i,k}, \va{u}_t^{i,k}\right)\right)_{t=0, \dots, T-1}\right), \\
			\intertext{where:}
			\va{\lambda}_{t+1}^k &= h_{\alpha_t^k}\left(\va{\lambda}_t^k, \va{\xi}_{t+1}\right), \qquad \forall t=0, \dots, T-1.
		\end{align*}
	\UNTIL{$\left\Vert \lambda^{k+1} - \lambda^k \right\Vert < \varepsilon$}
\end{algorithmic}
\end{algorithm}

\begin{rem}[Convexity of $\mathcal{S}$]
	The regression operator $\mathcal{R}^s$ is meant to be a sample-based approximation of the projection operator $\mathcal{R}$ on $\mathcal{S}$. Since the latter set can be non-convex, $\mathcal{R}\big(\va{\lambda}^{k+\frac{1}{2}}\big)$ is not necessarily unique: this may lead to numerical instabilities.
\end{rem}

\begin{rem}[Enhancement of $\mathcal{S}$]
It may be desirable to consider a larger set $\mathcal{S}$ in order to estimate more accurately the price process. For instance, one can extend relation \eqref{eqn:DualDynamics} in order to include more memory in the process:
\begin{equation*}
	\va{\lambda}_{t+1} = h_{\alpha_t}\left(\left(\va{\lambda}_\tau\right)_{\tau \leq t}, \left(\va{\xi}_\tau\right)_{\tau \leq t+1}\right),
\end{equation*}
However, this will in general increase the numerical complexity of DP in the solving of the subproblems.
\end{rem}

\begin{rem}[Another formulation]
	Alternatively, we could have considered a gradient algorithm that iterates directly on the parameters $\alpha$ of the dynamics \eqref{eqn:DualDynamics}. In this case, since we have no restrictions on $\alpha$, the feasible set would have been convex. Unfortunately, because the dynamics \eqref{eqn:DualDynamics} may be nonlinear with respect to $\alpha$, the dual function $\psi$ introduced in \S\ref{sec:PD} might be non-concave with respect to $\alpha$.
\end{rem}

\section{Numerical experiments} \label{sec:NumExp}

We tested this approach on a simple power management problem. On this small-scale example, we will be able to compare DADP results to those obtained by DP. Consider a power producer who owns two types of power plants:
\begin{itemize}
	\item Two hydraulic plants that are characterized at each time step $t$ by their water stock $\va{x}_t^i$ and power production $\va{u}_t^i$, and receive water inflows $\va{\xi}_{t+1}^i$, $i=1, 2$. These two units are subject to dynamic constraints but are cost-free;
	\item One thermal unit with a production cost that is quadratic with respect to its production $\va{u}_t^3$. There are no dynamics associated with this unit.
\end{itemize}
Using these plants, the power producer must supply a power demand $\va{d}_t$ at each time step $t$, over a discrete time horizon of $T=25$ time steps. All noises, i.e. the demand $\va{d}_t$ and the inflows $\va{\xi}_t^1$ and $\va{\xi}_t^2$ are chosen to be white noise processes.
\begin{figure}[ht]
\begin{tabular}{cc}
	\includegraphics[height=4.5cm]{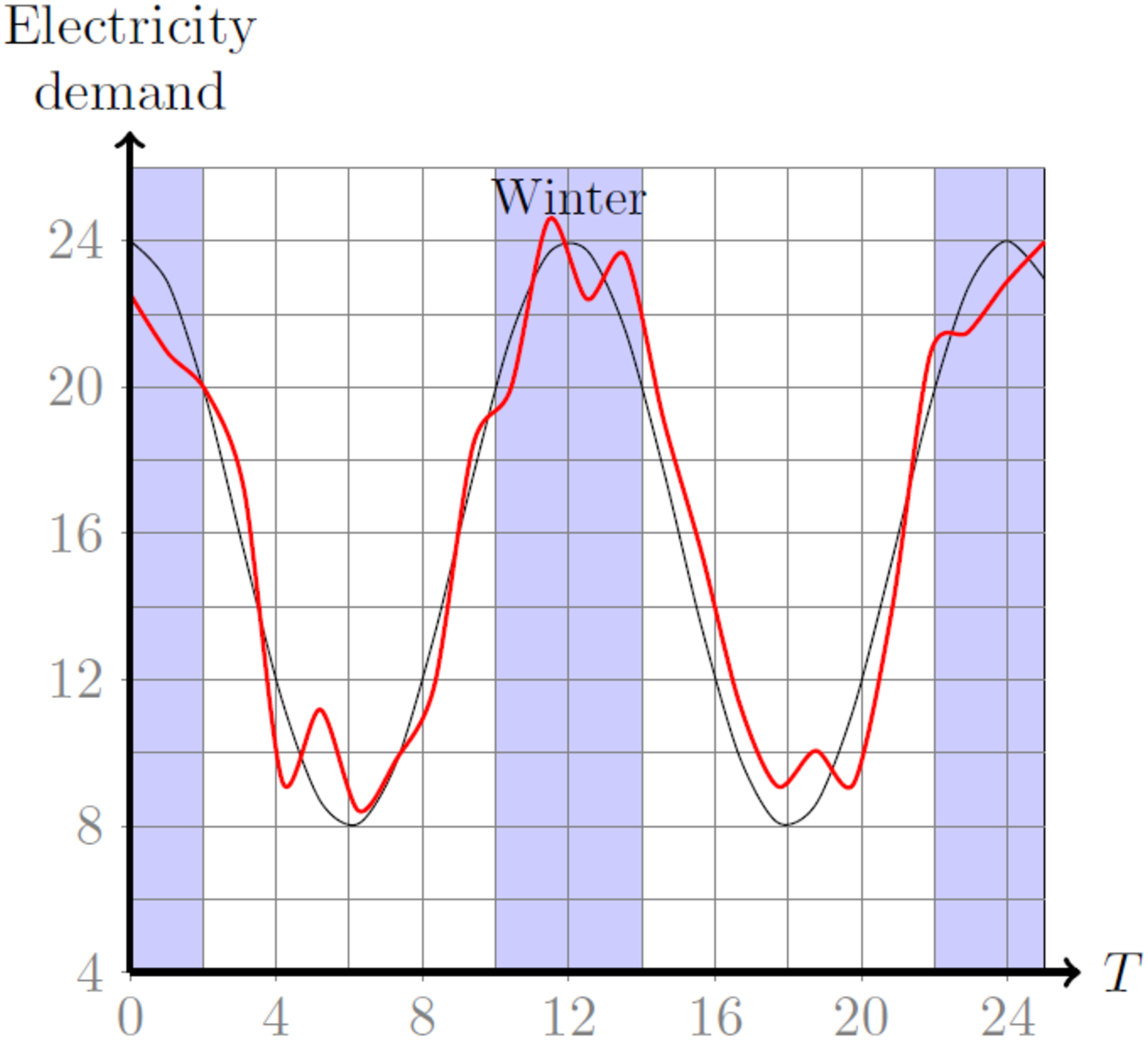} \\
	\includegraphics[height=4.5cm]{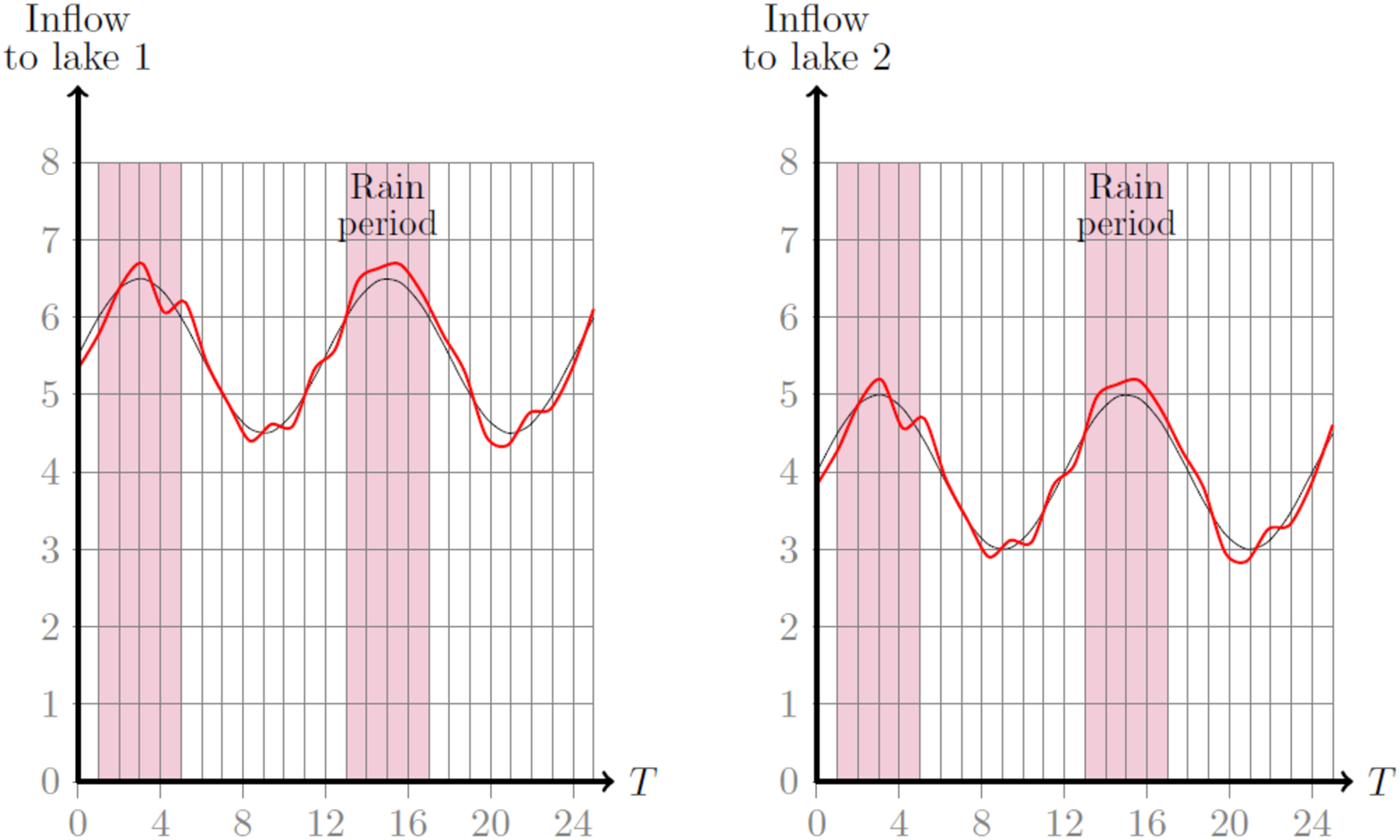}
\end{tabular}
\caption{\label{fig:DemandInflows}Mean of $\va{d}_t$, $\va{\xi}_t^1$ and $\va{\xi}_t^2$ over time (in red is one sample trajectory).}
\end{figure}

In the model we moreover impose small quadratic costs on the hydraulic power productions in order to ensure that, at least in the deterministic framework and without our approximation, the algorithm would build primal iterates that converge to the optimal solution of the original problem (see \cite{EkelandTemam92}). The problem reads\footnote{In this example, we consider two hydraulic plants with characteristics:
\begin{equation*}
	\underline{x}^1 = 0, \qquad \overline{x}^1 = 50, \qquad \overline{u}^1 = 6, \qquad K^1\left(x\right)=-7 x,
\end{equation*}
\begin{equation*}
	\underline{x}^2 = 0, \qquad \overline{x}^2 = 40, \qquad \overline{u}^2 = 6, \qquad K^2\left(x\right)=-12 x, \qquad \epsilon = 0.1
\end{equation*}
where $\underline{y}$ (resp. $\overline{y}$) denotes a lower (resp. upper) bound for variable $y$. Moreover, producing $u$ with the thermal plant costs $L_t\left(u\right) = u+u^2$.}:
\begin{subequations} \label{eqn:PTest}
\begin{align}
	\min_{\va{x}, \va{u}} \quad &\esper{\sum_{t=0}^{T-1} \left(\epsilon \left(\va{u}_t^1\right)^2 + \epsilon \left(\va{u}_t^2\right)^2 + L_t\left(\va{u}_t^3\right)\right) + K^1\left(\va{x}_T^1\right) + K^2\left(\va{x}_T^2\right)} \\
	\text{s.t.} \quad &\va{x}_{t+1}^i = \va{x}_t^i - \va{u}_t^i + \va{\xi}_{t+1}^i, \qquad \forall i=1,2, \quad \forall t=0,\dots, T-1, \\
	& \va{u}_t^1 + \va{u}_t^2 + \va{u}_t^3 = \va{d}_t, \qquad \forall t=0, \dots, T-1, \label{eqn:PTest3} \\
	&\underline{x}^i \leq \va{x}_t^i \leq \overline{x}^i, \qquad \forall i=1,2, \quad \forall t=1,\dots, T, \\
	&0 \leq \va{u}_t^i \leq \overline{u}^i, \qquad \forall i=1,2, \quad \forall t=0,\dots, T-1, \\
	&0 \leq \va{u}_t^3, \qquad \forall t=0, \dots, T-1, \\
	& \va{u}_t^i \text{ is } \sigma\big\{\va{d}_0, \va{\xi}_0^1, \va{\xi}_0^2, \dots, \va{d}_t, \va{\xi}_t^1, \va{\xi}_t^2\big\} \text{-measurable}, \quad \forall i=1,2,3.
\end{align}
\end{subequations}

In this problem, the state $\va{x}_t$ is two-dimensional, hence DP remains numerically tractable and we can use the DP solution as a reference. In order to use DADP, we choose an auto-regressive process for the Lagrange multipliers:
\begin{subequations} \label{eqn:ChosenShape}
\begin{align}
	\va{\lambda}_{t+1} &= \alpha_t \va{\lambda}_t + \beta_t \va{d}_{t+1} + \gamma_t, \\
	\va{\lambda}_0 &= \beta_0 \va{d}_0 + \gamma_0.
\end{align}
\end{subequations}

\begin{figure}[ht]
\begin{center}
	\includegraphics[width=0.8 \textwidth]{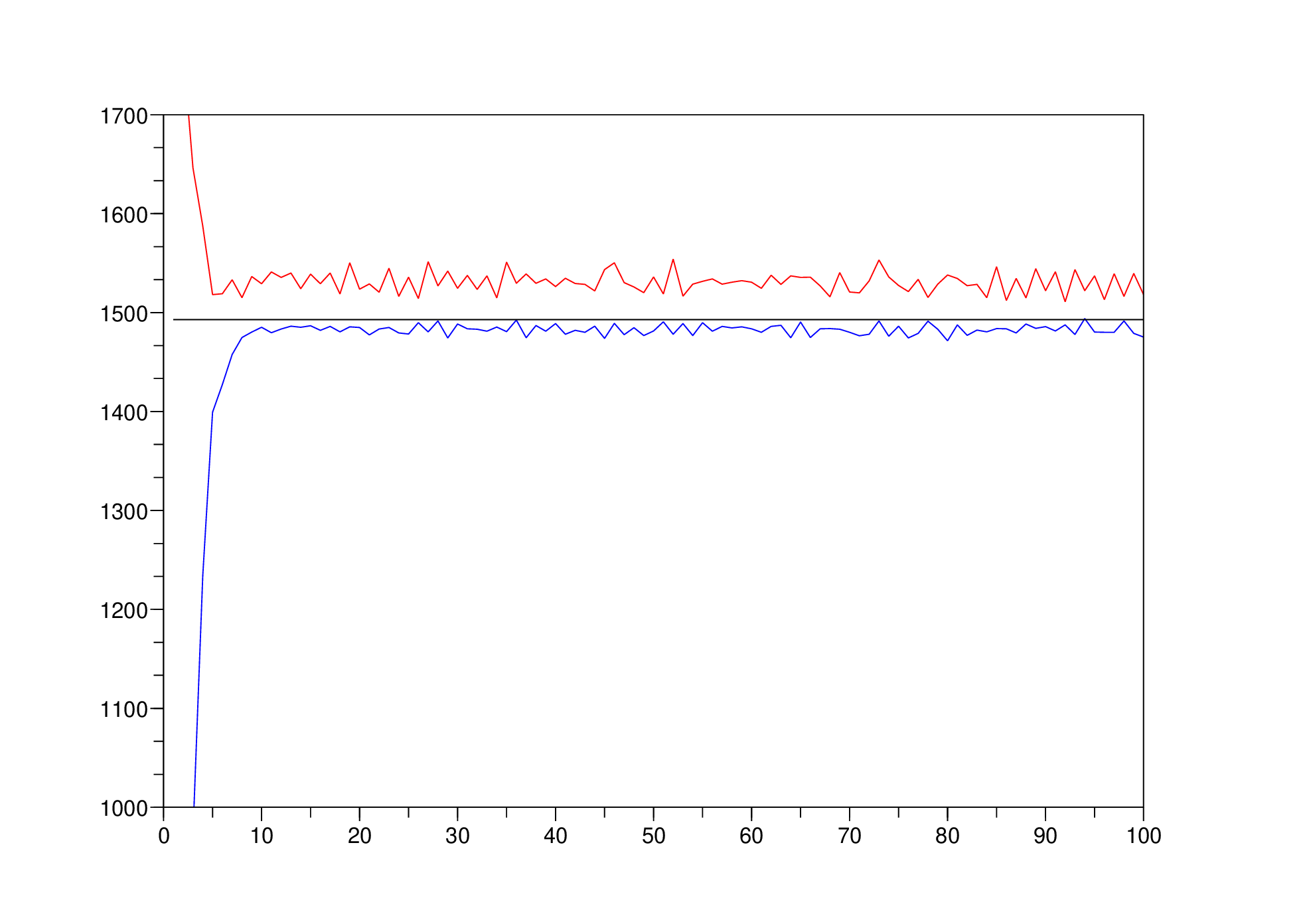}
\end{center}
\caption{\label{fig:DualValue}Value of the dual function (blue), of the primal function (red) and optimum (black) along with iterations.}
\end{figure}

We then perform the algorithm and depict its convergence in Figure \ref{fig:DualValue}. We first draw the values of the dual function $\psi$ introduced in \S\ref{sec:PD} along with iterations (lower curve) and observe that it converges to the optimal value of the original problem computed by DP.
Note that each value of $\psi$ is computed by Monte Carlo simulation over $10^3$ scenarios. We also draw the cost of the problem with all constraints satisfied (primal cost) at each iteration (upper curve). As explained in \S\ref{sec:Point}, DADP does not ensure that the coupling constraint \eqref{eqn:PTest3} is satisfied. To circumvent this difficulty, the thermal unit strategy is chosen in the simulation so as to ensure feasibility of the coupling constraint, i.e.:
\begin{equation} \label{eqn:SatisfyDemand}
	\va{u}_t^3 = \va{d}_t - \left(\va{u}_t^1+\va{u}_t^2\right).
\end{equation}
That is, DADP returns three strategies, for each of the hydraulic units and for the thermal unit. However, we use relation \eqref{eqn:SatisfyDemand} for the thermal strategy during simulations in order to ensure demand satisfaction.

Figure \ref{fig:DualValue} shows that the algorithm behaves well, in the sense that the value of the objective function converges quite quickly to a neighbourhood of the optimal value. However, even after 100 iterations the curve is still a bit noisy. This is because the price dynamics employed generates a non-convex-set of stochastic processes. Consequently, the least squares problem solved at each iteration is non-linear, and hence small variations in the actual gradient can result in large changes in the calculated gradient.

\begin{figure}[ht]
\begin{center}
\begin{tabular}{cc}
	 \includegraphics[width=5cm]{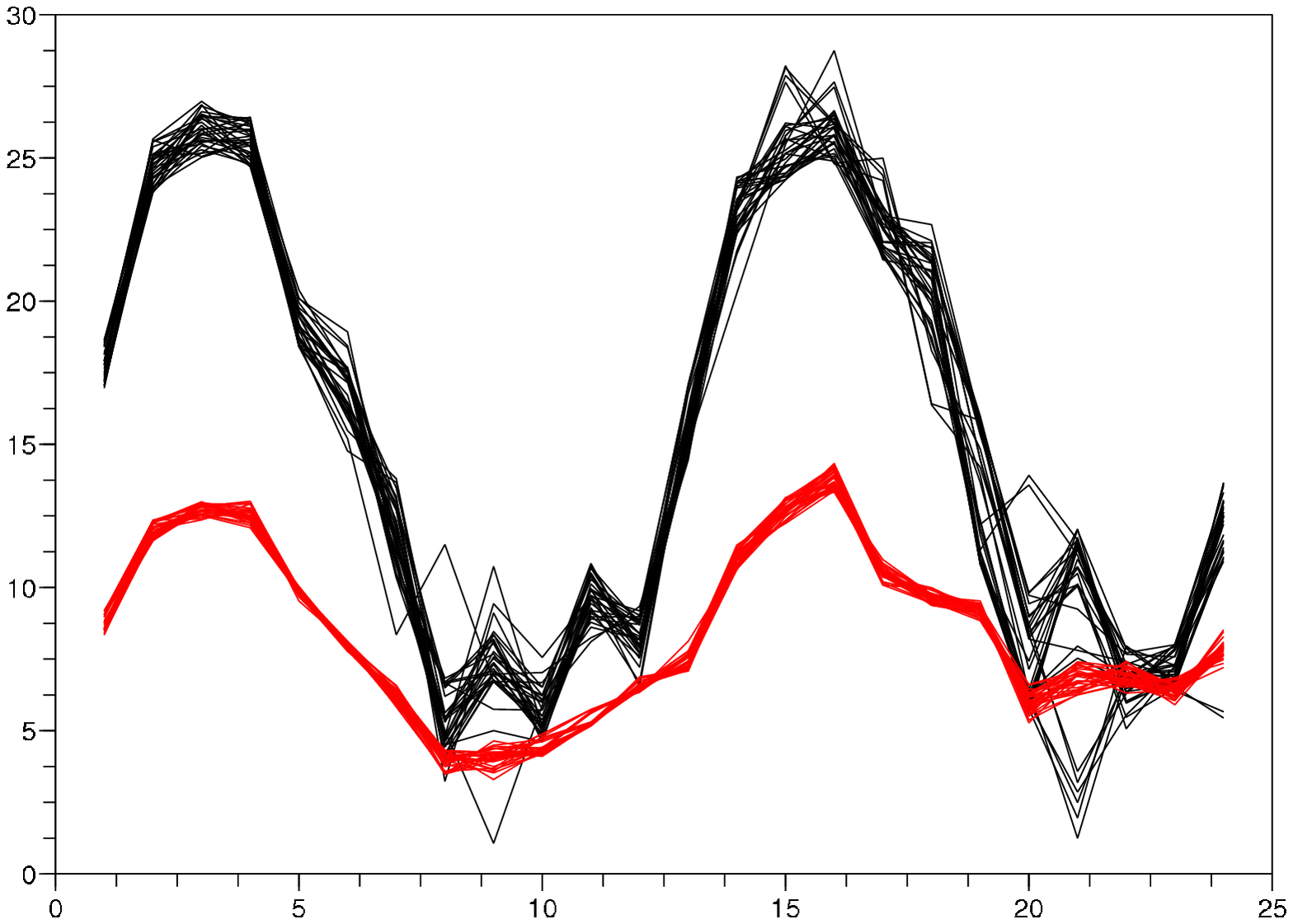}
	&\includegraphics[width=5cm]{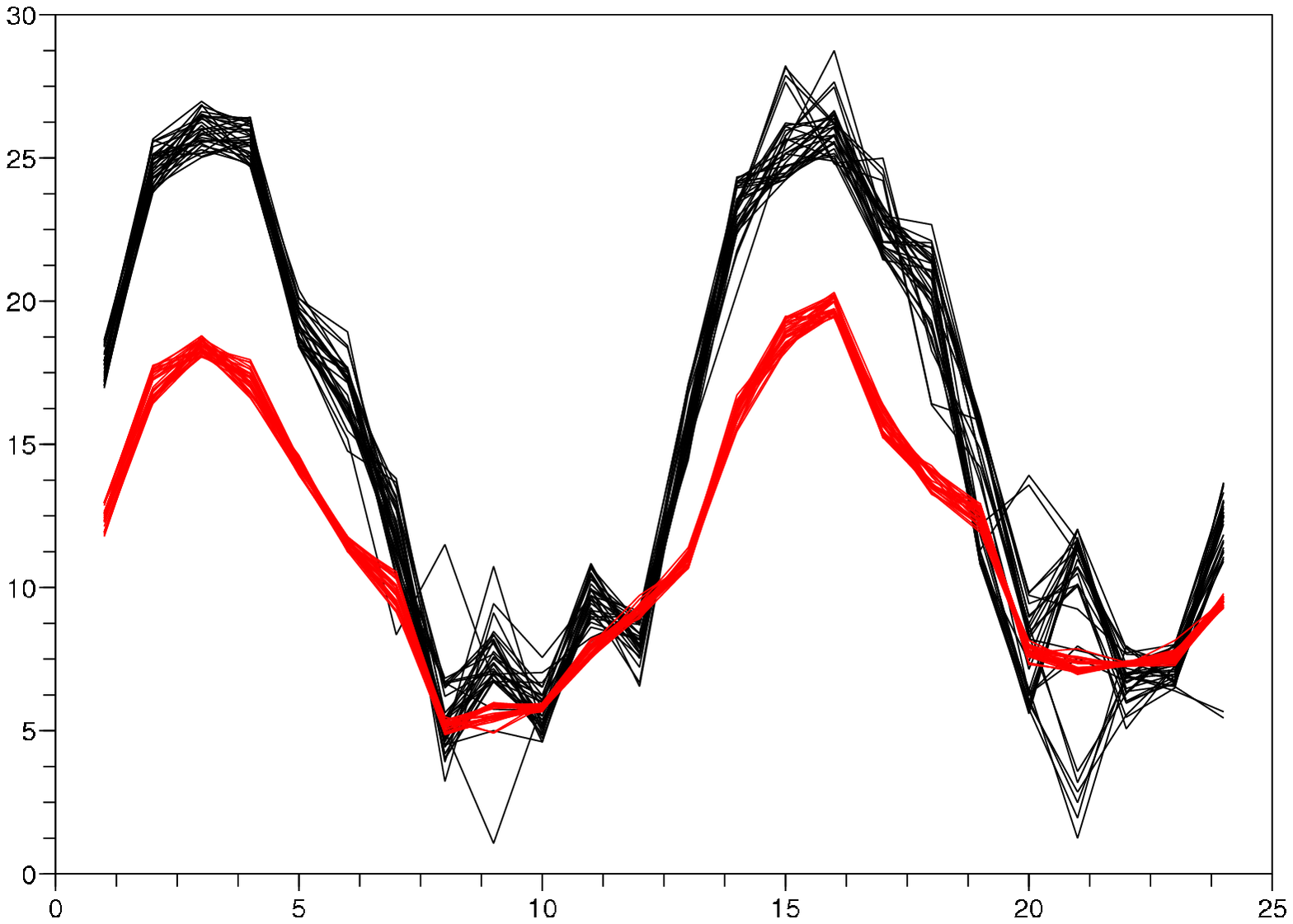} \\
	 \includegraphics[width=5cm]{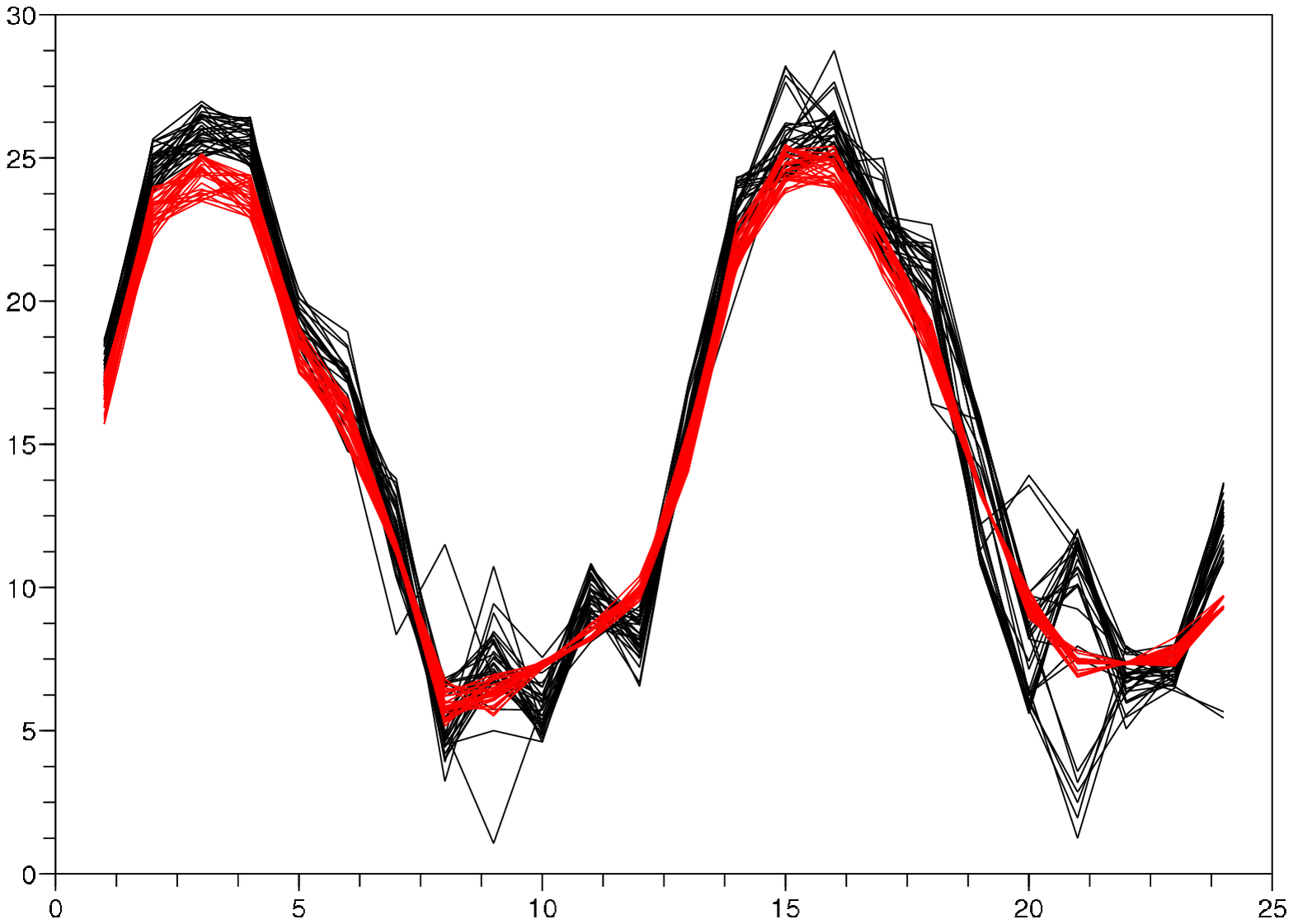}
	&\includegraphics[width=5cm]{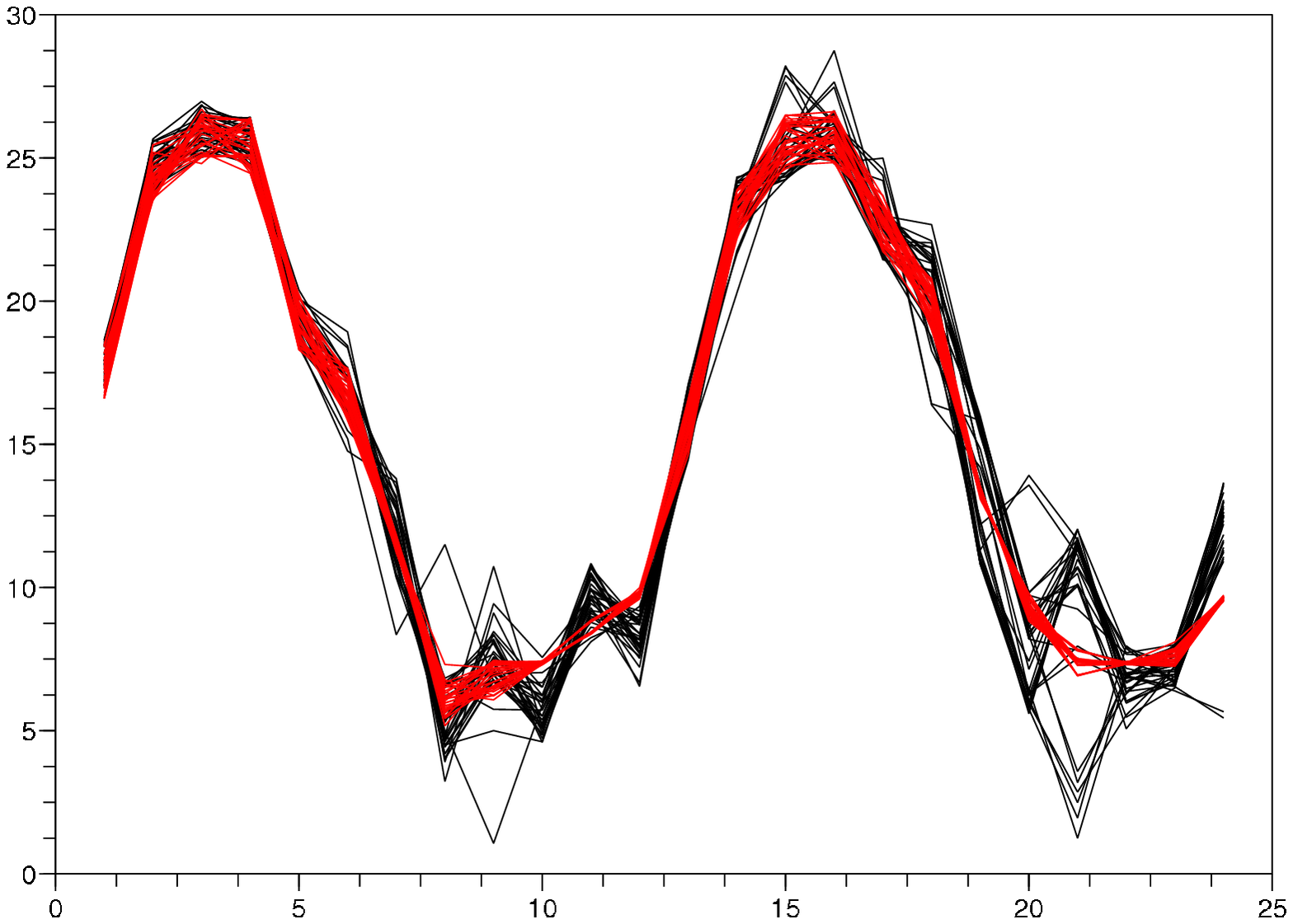}
\end{tabular}
\caption{\label{fig:PricesConvergence} Comparison, using 100 samples, of the approximate prices trajectories (red) with the optimal ones (black), after 10, 20, 50, and 90 iterations of the algorithm (from left to right and top to bottom).}
\end{center}
\end{figure}

The key to convergence in DADP is to obtain a dynamics for the Lagrange multipliers that accurately matches the optimal one. We have represented in Figure \ref{fig:PricesConvergence} the dynamics of the multipliers computed by DADP after 10, 20, 50 and 90 iterations, and those derived from DP. We observe that the approximate price dynamics issued from DADP satisfactorily converges to the optimal one.
This indicates that:
\begin{enumerate}
	\item the dual process converged, although the set of stochastic processes defined by \eqref{eqn:ChosenShape} is non-convex,
	\item there is no need to enhance the chosen dynamics \eqref{eqn:ChosenShape} for the multipliers in this particular problem.
\end{enumerate}

Note that the optimal prices derived from DP are obtained by numerically differentating the Bellman functions, hence the numerical instabilities we observe in the lower parts of the DP prices curves.

\begin{rem}[Numerical complexity]
We chose to validate the method on a two-dimensional power management problem, so we could compute a reference solution using DP. Note that there would be no additional difficulty in implementing the same algorithm with a larger number of hydraulic units, i.e. with a larger-dimensional state. The complexity of DADP grows linearly with the number of subsystems. The most time consuming part of the algorithm is solving the subproblems. However, this calculations can be easily parallelized on a computer, so that the time needed at each DADP iteration remains constant with respect to the number of subsystems.
\end{rem}

\section{Conclusion}

In this paper we present an approach, called Dual Approximate Dynamic Programming (DADP), to solve large-scale stochastic optimal control problems in a price decomposition framework, without discretizing randomness. This method employs classical duality theory results to solve decomposable large-scale systems, as it usually is in the deterministic framework. In order to be able to solve subproblems using DP, we suppose that the Lagrange multipliers obey some parameterized dynamics. The DADP algorithm then iterates on the parameters of these dynamics. What is original in this approach is the use of a dual variable in the optimal local feedback functions as an auxiliary variable that sums up the remaining part of the system.

On an example, we show that this approach is very attractive from a numerical point of view. Using rather simple dynamics for the multipliers, we obtained surprisingly good results with a small number of iterations. The main advantage of the method is that the complexity of the algorithm grows linearly with respect to the number of subsystems so that the curse of dimensionality is circumvented for the considered class of problems.

There are still several important theoretical questions. Since we constrain the dual variables to lie in some a priori chosen subset, we cannot state that the coupling constraints will be satisfied. Hence it would be useful to be able to evaluate the distance between the solution given by the heuristic and the feasible set; this would also give clues on how to choose well-suited dynamics for the dual variables on particular problems. Furthermore, the stochastic process subset on which we constrain the dual variables is possibly non-convex. In this context, it might be valuable to use more enhanced numerical methods for the update of the Lagrange multipliers. Further studies will be concerned with these issues.


\bibliographystyle{amsalpha}
\bibliography{biblio}

\end{document}